\magnification=\magstep1
\input amstex
\documentstyle{amsppt}

\define\defeq{\overset{\text{def}}\to=}
\define\ab{\operatorname{ab}}
\define\tp{\operatorname{tp}}
\define\Gal{\operatorname{Gal}}

\define\Inn{\operatorname{Inn}}
\define\Isom{\operatorname{Isom}}

\define\metab{\operatorname{metab}}
\define\Spec{\operatorname{Spec}}

\def\Aut{\operatorname{Aut}}
\def\Hom{\operatorname{Hom}}

\def\Out{\operatorname{Out}}

\NoRunningHeads
\NoBlackBoxes
\topmatter

\title
My view on and experience with Grothendieck's anabelian geometry
\endtitle

\author
Mohamed Sa\"\i di
\endauthor

\abstract This paper is based on the author's talk at the Grothendieck conference at Chapman university on 26-05-2022. It is not a survey of anabelian geometry but rather exposes some of the personal views and experiences of the author with the topic. 
\endabstract

\toc
\subhead
\S0. Introduction 
\endsubhead

\subhead
\S1. Following Galois
\endsubhead

\subhead
\S2. Anabelian geometry
\endsubhead

\subhead
\S 3. $m$-Step solvable anabelian geometry
\endsubhead

\subhead
\S 4. The Grothendieck philosophy
\endsubhead

\endtoc

\endtopmatter

\document
\subhead
\S 0. Introduction
\endsubhead
Alexander Grothendieck (1928-2014) was undoubtedly  one of the greatest mathematician of all times. The impact and depth of his work is still very hard to measure both by his contemporaries and successors alike.  When I was asked by the organisers of the Grothendieck conference to give a talk I was both humbled and a bit worried because of the greatness of Grothendieck's legacy. I also realised that I did spend the last 32 years of my active mathematical life (since I started my PhD) thinking, working, and pondering on mathematical topics and questions that interested Grothendieck and in which he himself had a substantial contribution, so I thought after all that I should be able to give the talk and give Grothendieck the credit he deserves.

\bigskip
For me Grothendieck was a mathematician well ahead of his time. After half a century since he stopped actively and publicly doing mathematics, a major research activity is ongoing in research topics that he initiated and launched. Researchers involved in these activities are still discovering and admiring the depth and impact of his ideas. 

\bigskip
The conference was actively attended by philosophers of science and mathematics and logicians, a fact I admired very much because of my naive interest in philosophy and its relationship with mathematics.
During the conference, listening and discussing with various people, I realised among others that the use by working mathematicians of the terms ``philosophy" and ``philosophical'' is rather casual and wouldn't correspond to how philosophers perceive these terms. I nonetheless ventured during my talk in the philosophy territory claiming that Grothendieck was a very special kind of mathematician, in comparison to the way other mathematicians do and perceive mathematics, and indeed that he was a``philosopher mathematician". 

\bigskip
I did also claim in my talk that Grothendieck, despite the substantial activity and the sheer number of mathematicians that surrounded him during his active public mathematical life, was nonetheless an isolated mathematician. Indeed it seems to me that some of the mathematical topics that were extremely high in Grothendieck's mathematical list of interests were not actively shared by his contemporaries and disciples. I very much think in particular about the theory of arithmetic fundamental groups and anabelian geometry. 

\bigskip
No living mathematician can probably fairly and accurately speak on the whole of the mathematical legacy of Grothendieck. In my talk I chose to speak on the part of this legacy that I certainly know best which is his outstanding contribution to Galois theory and his launch of anabelian geometry.  

\bigskip
One philosophical thread in my talk was the following general, both mathematical and philosophical (in the sense used by mathematicians), question:

\bigskip
{\it What is the precise gap (say in Galois theory broadly speaking) between commutative and non-commutative mathematics, as well as between linear and nonlinear mathematics}?

\subhead
\S1. Following Galois
\endsubhead
In his outstanding work, \'Evariste Galois (1811-1832) solved the ancient famous problem of solvability of polynomial equations by radicals. 

\bigskip
As one learns in a beginners Galois theory course, Galois proved that a polynomial in one variable with integer coefficients is solvable by radicals if and only if the Galois group of this polynomial is a solvable (finite) group.
This result, proved two centuries ago, seems to me to be the first theorem in what one might call "modern mathematics",  whereby one proves a mathematical statement $A$ regarding a mathematical object  
$\tilde A$, by relating object $\tilde A$ to a mathematical object $\tilde B$ of a completely different nature, and establishing an equivalence between statement $A$ and a mathematical statement $B$ related to object $\tilde B$. In the case of Galois, $\tilde A$ is a polynomial in one variable $f(X)$, and $\tilde B$ is the finite group $\Gal(f(X))$: the Galois group of $f(X)$, statement $A$ is the solvability of $f(X)$ by radicals, and statement $B$ is $\Gal(f(X))$ is a solvable group.

\bigskip
Almost a century after Galois' work the notion of the absolute Galois group of a field has been introduced together with its profinite topology. At that time the inverse Galois problem, which asks
what finite groups occur as Galois groups of polynomials with rational coefficients, had already been initiated by Hilbert and his collaborators.

\bigskip
The problem of explicitly determining the structure of the Galois group $G_{K}$ of a field $K$, say by providing explicit generators and relations of this Galois group as a profinite group, turned out to be an extremely difficult problem in general. Even in the ``simplest case" of the prime field $\Bbb Q$ of rational numbers, the explicit determination of the structure of the Galois group $G_{\Bbb Q}$ seems to be out of reach for the time being. 
The situation for a finitely generated field over $\Bbb Q$ (or the prime field $\Bbb F_p$ of characteristic $p>0$, apart from finite fields) does not fare better.

\bigskip
In the case of number fields, i.e. finite field extensions of $\Bbb Q$, various approaches have been adopted by mathematicians in order to understand Galois groups, circumventing 
the problem of explicitly determining their structures. The most successful approach to date is class field theory (CFT), which occupied number theorists in the first half of last century.
This theory explicitly determines the structure of the maximal abelian quotient $G_K^{\ab}$ of the Galois group $G_K$ of a number field in terms of some data and structures which depend solely on $K$. This theory has been the precursor of the famous Langlands programme which is a central research programme in number theory and arithmetic geometry generating substantial ongoing research activities. Another extremely successful approach is the study of Galois representations of Galois groups of number fields initiated by Weil, Shimura, Grothendieck, Serre, and Deligne and which culminated in the works of Faltings and Wiles on the Mordell and Fermat conjectures. These (linear) Galois representations 
provided by  \'etale cohomology are not faithful; this is why Grothendieck was so struck by Belyi's result on the Galois representation on the geometric fundamental group of the projective line minus three points (with its rather short elementary proof!) which shows that there are faithful ``nonlinear Galois representations" of $G_{\Bbb Q}$.

\bigskip
Grothendieck entered the scene of Galois theory during the second half of the 1950s decade. The first volume SGA1 of his series of S\'eminaires de G\'eom\'etrie alg\'ebrique , which appeared in print in 1960, is devoted to the definition and study of basic properties of \'etale fundamental groups. 

\bigskip
In [Grothendieck] Grothendieck generalised the notion of absolute Galois group of a field  to schemes.
For a connected scheme $X$, and a fixed geometric base point $x$ of $X$, Grothendieck defined the Galois category of pointed finite \'etale covers of $X$, and the \'etale fundamental group $\pi_1(X,x)$ of $X$ with base point $x$ as the automorphism group of a certain functor on this Galois category. Serre and Grothendieck are credited for defining the right notion of cover in algebraic geometry which is the notion of finite \'etale cover, whose properties are studied in detail in [Grothendieck].

\bigskip
One of the first achievements of Grothendieck in algebraic geometry is the description of the structure of the prime-to-characteristic geometric \'etale fundamental group of an algebraic curve in positive characteristic. 

\bigskip
Let $X$ be a projective smooth and connected algebraic curve of genus $g\ge 0$ over an algebraically closed field $k$ with geometric point $x$. In the case $k=\Bbb C$ is the field of complex numbers,
and as an application of the Riemann existence theorem, one shows that the \'etale fundamental group $\pi_1(X,x)$ is isomorphic to the profinite completion
$\Gamma_g$ of the topological fundamental group $\pi_1(X(\Bbb C))^{\tp}$ of the topological space underlying the Riemann surface $X(\Bbb C)$ associated to $X$. The latter is isomorphic to a surface group of genus $g$ and is generated by $2g$ generators $\{a_1,\ldots,a_g,b_1,\ldots,b_g\}$ subject to the unique relation $\prod a_ib_ia_i^{-1}b_i^{-1}=1$.
In case $k$ has characteristic zero a Lefschetz principle argument implies that a similar result holds. In case $\text {char}(k)=p>0$, Grothendieck proved in [Grothendieck] that the maximal prime-to-$p$ quotient $\pi_1(X,x)^{(p')}$ of  $\pi_1(X,x)$
is isomorphic to the maximal prime-to-$p$ quotient $\Gamma_g ^{(p')}$ of $\Gamma_g$. In particular, the structure of
$\pi_1(X,x)^{(p')}$ as a profinite group depends only on the genus $g$ of $X$. Grothendieck proved this result as a consequence of his theory of specialisation for fundamental groups
and his theorem of existence in formal geometry.

\bigskip
During the academic year 1994-1995 I met the Dutch mathematician Jacob Pieter Murre at the University of M\"unster in Germany, where I was a post-doc and he was visiting Christopher Deninger. He informed me that the above result proven by Grothendieck was conjectured by Weil. He further told
me that while he was a fellow at Chicago University during the years 1954-1956, and despite the frosty relationship between Weil and Grothendieck (we learnt a great deal about this relationship during various talks at the conference, which it seems goes back to the days when both were members of the Bourbaki group), Weil urged Murre to go and study with Grothendieck the new algebraic geometry developed by him, arguing that the new language of algebraic geometry must be very powerful since it allowed Grothendieck to prove the result mentioned above.

\bigskip
In the case where $X\to \Spec K$ is a geometrically connected and smooth algebraic variety over a field $K$, $x$ is a geometric point of $X$ with values in its generic point, $\overline K$ is the algebraic closure of $K$ determined by $x$, $\overline X=X\times _K{\overline K}$ is the geometric fibre of $X$, the fundamental group $\pi_1(X,x)$ sits naturally in the following exact sequence of profinite groups:
$$1\to \pi_1(\overline X,x)\to \pi_1(X,x)\to G_K\to 1\tag *$$ 
where $G_K=\Gal(\overline K/K)$ is the Galois group of $\overline K/K$. 

\bigskip
The above exact sequence (*) induces naturally an outer Galois representation
$$\rho_X:G_K\to \Out (\pi_1(\overline X,x)),\tag **$$
where 
$$\Out (\pi_1(\overline X,x))=\Aut (\pi_1(\overline X,x))/\Inn(\pi_1(\overline X,x))$$ 
is the outer automorphism group of the profinite group $\pi_1(\overline X,x)$.

\bigskip
In the case where the group $ \pi_1(\overline X,x)$ is centre-free one can reconstruct the exact sequence (*) from the outer representations (**) and the exact sequence
$$1\to  \pi_1(\overline X,x)\to \Aut  \pi_1(\overline X,x) \to \Out  (\pi_1(\overline X,x))\to 1$$
by pull-back via $\rho_X$.

\bigskip
An important example is the case $K=\Bbb Q$ and $X=E$ is an elliptic curve. In this case  $\pi_1(\overline X,x)$ is the Tate module of $E$ and the (linear) Galois representation $\rho_X$ has been extensively investigated since the work of Serre in the 1960's and 1970's.

\bigskip
Two further important examples are the cases where $K=\Bbb Q$, $X=E\setminus \{0_E\}$ is a once punctured elliptic curve, and $X=\Bbb P^1\setminus\{0,1,\infty\}$. In both cases 
$\pi_1(\overline X,x)$ is profinite free on two generators, though the Galois representations $\rho_{E\setminus \{0_E\}}$ and $\rho _{\Bbb P^1\setminus \{0,1,\infty\}}$ are not isomorphic!

\bigskip
The study of the above exact sequence (*) as well as the outer Galois representation (**) sat extremely high in the list of mathematical interests of Grothendieck. Luc Illusie, a former student 
of Grothendieck, informed me on several occasions that Grothendieck told him in the year 1966 he had ideas how to solve Fermat's conjecture using the outer Galois  
representation $\rho _{\Bbb P_1\setminus \{0,1,\infty\}}$!

\subhead
\S2. Anabelian Geometry
\endsubhead
The terms anabelian and anabelian geometry appeared for the first time in the famous letter that Grothendieck sent to Faltings in the early 1980's after the latter proved Mordell conjecture. In this letter Grothendieck formulated a couple of conjectures: the birational anabelian conjecture for finitely generated fields, the anabelian conjecture for hyperbolic curves over finitely generated fields, and the section conjecture (I will not discuss the section conjecture in this paper). 

\bigskip
Let $K$ be a characteristic zero finitely generated field. The birational anabelian conjecture for finitely generated fields and the anabelian conjecture for hyperbolic curves over finitely generated fields read as follows.

\bigskip
$\bullet$ {\bf AN1} ({\bf The birational anabelian conjecture for finitely generated fields})
Let $F$ and $L$ be finitely generated fields over $K$. The natural set-theoretic map
$$\Hom _K(F,L)\to \Hom _{G_K}(G_L,G_F)/ \Inn G_F,$$
where $\Hom _{G_K}(G_L,G_F)/\Inn G_F$ is the set of continuous homomorphisms $G_L\to G_F$ which 
are compatible with the projections onto $G_K$ modulo inner automorphisms of $G_F$, is a {\bf bijection}.

\bigskip
$\bullet$ {\bf AN2} ({\bf The anabelian conjecture for hyperbolic curves over finitely generated fields})
Let $X$ and $Y$ be hyperbolic curves over $K$. The natural set-theoretic map
$$\Hom _K(X,Y)\to \Hom _{G_K}(\pi_1(X,\star),\pi_1(Y,\star))/ \Inn \pi_1(Y,\star),$$
where  $\Hom _{G_K}(\pi_1(X,\star),\pi_1(Y,\star))$
is the set of continuous homomorphisms $\pi_1(X,\star)\to \pi_1(Y,\star)$ which 
are compatible with the projections onto $G_K$ modulo inner automorphisms of $\pi_1(Y,\star)$; $(\star)$ are base points, is a {\bf bijection}.

\bigskip
The essence of these conjectures is that one can reconstruct embeddings between finitely generated fields, as well as morphisms between hyperbolic curves, from profinite group homomorphisms between absolute Galois groups and \'etale fundamental groups.

\bigskip
It is not very clear to me why Grothendieck wrote his letter to Faltings and shared with him his conjectures. One such reason, which was mentioned by 
Grothendieck in the letter, is the similarity of these conjectures with the Tate conjecture on Galois representations on Tate modules of abelian varieties 
which Faltings proved in his way towards proving Mordell. The philosophy underlying both conjectures can be very roughly summarised as follows:

\bigskip
{\it Instances such that the Galois representation $\rho_X$, whereby both arithmetic and topology mix with each other, should lead to extremely rigid situations!}

\bigskip
Another possibility might be that Grothendieck hoped Faltings would work on his conjectures and make progress in anabelian geometry. 

\bigskip
The precision and care taken by Grothendieck to formulate his conjectures is astonishing: the conjectures turn out to be true exactly in the form in which they were formulated!
It is also extraordinary that Grothendieck formulated the conjectures without providing any evidence in the form of an example or alike. In fact, even now, it is very hard to produce examples to support most statements in anabelian geometry. It seems to me Grothendieck did spend ample time (possibly many years) working on these conjectures 
and their formulation.
It is not clear also when exactly Grothendieck started thinking about these conjectures.

\bigskip
The isom-form of the birational anabelian conjecture, whereby one replaces field embeddings by field isomorphisms and homomorphisms between Galois groups by isomorphisms,
reads as follows (same assumptions as in {\bf AN1}): the map 
$$\Isom _K(F,L)\to \Isom _{G_K}(G_F,G_L)/ \Inn G_L $$
is a bijection. 

\bigskip
Similarly the  isom-form of the anabelian conjecture for hyperbolic curves
reads as follows (same assumptions as in {\bf AN2}): the map 
$$\Isom _K(X,Y)\to \Isom _{G_K}(\pi_1(X,\star),\pi_1(Y,\star))/ \Inn \pi_1(Y,\star)$$
is a bijection.

\bigskip
It turns out that the isom-form of the birational anabelian conjecture for number fields, i.e. finite extensions of the prime field $\Bbb Q$, was proven by Neukirch and Uchida a decade before Grothendieck wrote his letter to Faltings.
This fact is rather impressive because it shows evidence for the birational anabelian conjecture as formulated by Grothendieck. Moreover Neukirch and Uchida laid in their work the technical foundations which led to later breakthroughs in anabelian geometry. It is very likely Grothendieck was not aware of the work of Neukirch and Uchida for otherwise he would certainly have mentioned it in his letter to Faltings (there is no apparent reason why he would not have done so). 

\bigskip
Breakthrough in anabelian geometry came in the 1990s. The isom-form of the birational anabelian conjecture was proven by Pop (and Spiess in the case of Kronecker dimension $2$) building on the ideas of Neukirch and Uchida. The isom-form of the anabelian conjecture for hyperbolic curves was proven by Nakamura, Tamagawa, and Mochizuki. Subsequently Mochizuki proved a generalised version of the conjectures {\bf (AN1)} and {\bf (AN2)}, whereby one replaces $K$ by any sub-$p$-adic field (a subfield of a finitely generated field extension of a $p$-adic local field), using techniques from $p$-adic Hodge theory.

\bigskip
One consequence of the Grothendieck anabelian conjectures is that suitable natural functors
$$\Gal : \{\text {Finitely Generated Fields}\} \to \{\text {Profinite Groups}\}$$ 
$$\ \ \ \ L \ \mapsto \ G_L$$
and

$$\pi_1 : \{\text {Hyperbolic Curves}\} \to \{\text {Profinite Groups}\}$$ 
$$\ \ \ \ X\ \mapsto \ \pi_1(X,*)$$
are fully faithful. Here  $\{\text {Finitely Generated Fields}\}$,  $\{\text {Hyperbolic Curves}\}$, and $\{\text {Profinite Groups}\}$ denote the categories of 
 finitely generated fields,  hyperbolic curves, and profinite Groups respectively.

\bigskip
As was pointed out by Grothendieck in his letter to Faltings not only it is important to prove faithfulness of these functors but it is also important to determine the essential image of these functors. In other words what is the structure of Galois groups of finitely generated fields and \'etale fundamental groups of hyperbolic curves over finitely generated fields?
As mentioned earlier the task of determining these structures seems out of reach for the time being.

\bigskip
I personally started thinking about anabelian geometry some $25$ years ago and was always puzzled by the following question:

\bigskip
{\it How can a number theorist or an arithmetic geometer interested in concrete computations and applications use the main results of anabelian geometry?}

\bigskip
I personally am not aware of any concrete application of a computational nature of the Grothendieck anabelian conjecture. The reason, I already mentioned, is that we do not know the precise structure of Galois groups of finitely generated fields and \'etale fundamental groups of hyperbolic curves over finitely generated fields. The main issue I did consider was the following: 

\bigskip
{\it How to improve this situation?} 

\bigskip
One related issue is

\bigskip
{\it What is the precise meaning of the term anabelian?} 

\bigskip
Grothendieck carefully refrained from giving a precise definition of this term. Anabelian algebraic varieties
can be defined roughly as those varieties for which a statement as in {\bf AN2} holds. Thus hyperbolic curves over finitely generated field are anabelian varieties. It is still not clear what are precisely those algebraic varieties of dimension $>1$ which are anabelian.

\bigskip
The following is an extract of the paper [Mochizuki-Nakamura-Tamagawa] suggesting a ``definition" of the term anabelian:

\bigskip
{\it Here, the term ``anabelian algebraic variety" means roughly ``an algebraic variety whose geometry is controlled by its fundamental group, which is assumed to be 
{\bf far from abelian} "}.

\bigskip
This quotation suggests the term anabelian is perceived by experts as meaning ``far from being abelian". Anabelian geometry is thus the study of fields and algebraic varieties whose Galois groups and arithmetic fundamental groups respectively are ``far from being abelian", and the extent to which these fields and algebraic varieties can be reconstructed from their Galois groups and \'etale fundamental groups respectively.

\bigskip
This prompts the following question:

\bigskip
{\it What is the precise meaning for a profinite group to be {\bf far from abelian}?}

\bigskip
Grothendieck, as far as I know, never used the expression ``far from abelian". However he stressed the condition of centre-freeness for a profinite group which was generalised by Mochizuki who introduced the notion of slimness. A profinite group is slim if each of its open subgroups is centre-free.
Some experts argue and use the slimness condition as an adequate meaning for ``being far from abelian". I personally was never completely satisfied of this interpretation. 

\bigskip
In the last few years, together with my collaborator Akio Tamagawa, we discovered a new version of anabelian geometry which rather surprisingly suggests the intuition behind the above interpretation of the meaning of anabelian is rather a {\bf false} intuition.

\bigskip
\subhead
\S 3. $m$-Step solvable anabelian geometry
\endsubhead
For a profinite group $G$ define the profinite derived series 
$$...\subseteq G(i+1)\subseteq G(i)\subseteq...\subseteq G(1)\subseteq G(0)=G$$
where 
$$G(i+1)=\overline {[G(i),G(i)]},$$
for $i\ge 0$, is the closure of the $(i+1)$-th commutator subgroup $[G(i):G(i)]$.

\bigskip
The quotient
$$G^i\defeq G/G(i)$$
is the $i$-th {\bf step solvable quotient} of $G$. Thus $G^1=G^{\ab}$ is the maximal abelian quotient of $G$, $G^2=G^{\metab}$ is the maximal metabelian quotient of $G$, $\ldots$. 

\bigskip
For $j>i$ we have the following commutative diagram of exact sequences:
$$
\CD
1@>>> G(i)@>>> G @>>> G^i @>>>1\\
@. @VVV @VVV @V{\text id}VV\\
1@>>> G[j,i]@>>> G^j @>>> G^i @>>>1\\
\endCD
$$
where $G[j,i]$ is defined so that the lower sequence is exact.

\bigskip
Note that if $K/\Bbb Q$ is a finite field extension, and $m\ge1$, then the structure of $G_K^m$ can be in principle approached by class field theory.
In Iwasawa theory for example one describes explicitly the structure and arithmetic of certain quotients of $G_K^2$.

\bigskip
During the summer of 2017, together with Tamagawa, we proved the following result ([Sa\"\i di-Tamagawa]).

\proclaim {Theorem A} (Sa\"\i di-Tamagawa\ 2017) Let $L$ and $F$ be number fields. Then the natural map 
$$\Isom (F,L)\to \Isom (G_L^m,G_F^m)/\Inn G_F^m$$
is a bijection for all $m\ge 3$.
\endproclaim

\bigskip
Later we proved a similar result for finitely generated fields $L$ and $F$ with the same Kronecker dimension $d$ whereby $3$ is replaced by a ``small" constant depending on $d$.
A similar result is expected for hyperbolic curves over finitely generated fields whereby $\pi_1(X)$ and $\pi_1(Y)$ in {\bf AN2} are replaced by $\pi_1(X)^m$ and $\pi_1(Y)^m$,
respectively, for a small value of $m$.

\bigskip
Theorem A is rather surprising in many respects. 

\bigskip
First, it suggests that as far as anabelian geometry is concerned one doesn't need to know the structure of full absolute Galois groups of finitely generated fields but rather the structure of the quotients $G_K^m$ for ``small" values of $m$, which as mentioned above can be approached via class field theory.

\bigskip
Theorem A {\bf reconciles} anabelian geometry with class field theory which is a rather explicit and amenable to computation theory. This fact seems to be against the intuition of Grothendieck who, as we learnt during the conference, had a rather cold experience with class field theory.

\bigskip
Theorem A suggests that the {\bf anabelian} world is rather close to the {\bf abelian} world! 
This rather surprising fact, which suggests that in anabelian geometry the gap between the non-commutative and the commutative is not so large, goes against the intuition of most 
anabelian geometers, number theorists, and arithmetic geometers alike.

\bigskip
The idea that ``truncated" versions of Galois groups should be enough to encode the geometry of function fields is due to Bogomolov, who some 35 years ago conjectured (and provided important hints of proof) that two-step nilpotent Galois groups should be enough to reconstruct function fields of transcendence 
degree $\ge2$ over algebraically closed fields. Theorem A is rather arithmetic in nature. In fact it has been recently shown that two-step nilpotent Galois groups do not 
encode the isomorphy type of number fields (see [Koymans-Pagano]). The fact that an arithmetic truncated version of anabelian geometry as in Theorem A holds and the relation to class field theory were a big surprise and discovery to the author!

\bigskip
\subhead
\S 4. The Grothendieck philosophy
\endsubhead
The history of Theorem A is rather interesting. 

\bigskip
One day of May 2017, while in my office browsing on arXiv, I glanced into a recently posted paper by
Cornelissen, de Smit, Li, Marcolli, and Smit: ``Characterization of global fields by Dirichlet L-series". Immediately after reading the introduction, and in a flash instant, I felt very strongly that a statement as in Theorem A should be true. A couple of weeks later I visited Akio Tamagawa in Kyoto and shared with him my feeling. After 2 months, during which we met once a week for  mathematical discussions, we had proved Theorem A.

\bigskip
During the lockdowns of 2020 and 2021 due to Covid I started pondering on the following question:

\bigskip
{\it What exactly is the process that occurred in my mind and which led to the flash instant of May 2017 leading to the discovery of the statement of Theorem A?}

\bigskip
Indeed I didn't have any evidence in the form of an example or mathematical instance which could lead to Theorem A. More generally:

\bigskip
{\it How do mathematicians discover statements of theorems without any evidence in the form of examples or computations or alike?}

\bigskip
During the preparation of my talk at the Grothendieck conference, and while reading a paragraph in Grothendieck's letter to Faltings, 
I realised exactly what mental process went in my mind and which led to the flash instant of May 2017. 

\bigskip
What I realised is that in this mental process I had embraced  ``Grothendieck philosophy". 

\bigskip
I will not venture and  define what the {\bf mathematical philosophy of Grothendieck} is precisely but I will rather invite the reader to read and ponder on the following paragraph
from Grothendieck's letter to Faltings which I mentioned above:

\bigskip
{\it ``What my experience of mathematical work has taught me again and again, is that the {\bf proof always springs from the insight}, and not the other way round ? and that the insight itself has its source, first and foremost, in a {\bf delicate and obstinate feeling of the relevant entities and concepts and their mutual relations}. The guiding thread is the inner coherence of {\bf the image which gradually emerges from the mist}, as well as its consonance with what is known or foreshadowed from {\bf other sources} - and it guides all the more surely as the ``exigence" of coherence is stronger and more delicate"}.

\bigskip
{\it ``the image which gradually emerges from the mist"}, that's exactly what happened to me on that day of May 2017!

\bigskip
Grothendieck has left this world. He is, and will be, hugely missed by the mathematical community. Historians and mathematicians will still be debating his legacy for long time to come.
The question however for me is:

\bigskip
{\bf How can we benefit more from Grothendieck Today?}

\bigskip
My take on this is that we should:

\bigskip
\item {$\bullet$} Embrace more his mathematical philosophy in our way of doing research in mathematics, and reconcile his philosophy with ``practical mathematics".

\bigskip
\item {$\bullet$}  Embrace more his mathematical philosophy in our way of teaching mathematics.

\bigskip
$$\text{References.}$$

\smallskip
\noindent
[Grothendieck] Grothendieck, A., Rev\^etements \'etales et groupe fondamental, Lecture 
Notes in Math. 224, Springer, Heidelberg, 1971.

\smallskip
\noindent
[Koymans-Pagano] Koymans, P., Pagano, C., Tow step nilpotent extensions are not anabelian, arXiv:2301.10342, 2023.

\smallskip
\noindent
[Mochizuki-Nakamura-Tamagawa] Mochizuki, S., Nakamura, H., and Tamagawa, A., The Grothendieck Conjecture on the Fundamental Groups of Algebraic Curves,
Sugaku Expositions, Volume 14, Number 1, 31-53, June 2001.

\smallskip
\noindent
[Sa\"\i di-Tamagawa] Sa\"\i di, M., Tamagawa, A., The $m$-step solvable anabelian geometry of number fields, J. Reine Angew. Math. (Crelle), DOI: 10.1515/crelle-2022-0025.

\bigskip
\noindent
Mohamed Sa\"\i di

\noindent
College of Engineering, Mathematics, and Physical Sciences

\noindent
University of Exeter

\noindent
Harrison Building

\noindent
North Park Road

\noindent
EXETER EX4 4QF

\noindent
United Kingdom

\noindent
M.Saidi\@exeter.ac.uk

\end
\enddocument